\documentclass[preprint,12pt]{elsarticle}

\usepackage{bbm}
\usepackage{amsmath,amssymb,amsfonts,amsthm,
latexsym, amscd, amsfonts, epsfig, eepic,epic}
\usepackage{mathrsfs}
\usepackage{color}
\usepackage{eucal}%    caligraphic-euler fonts: \mathcal{ }
\usepackage{eufrak}%   frak-euler        fonts: \mathfrak{ }
\usepackage[all]{xypic}
\usepackage{xspace}

\theoremstyle{plain}
\newtheorem{theorem}{Theorem}
\newtheorem{algorithm}{Algorithm}
\newtheorem{corollary}{Corollary}

\newtheorem{lemma}{Lemma}

\theoremstyle{definition}

\theoremstyle{example}

\theoremstyle{remark}
\numberwithin{equation}{section}
\begin{document}

\begin{frontmatter}

%% Title, authors and addresses

%% use the tnoteref command within \title for footnotes;
%% use the tnotetext command for theassociated footnote;
%% use the fnref command within \author or \address for footnotes;
%% use the fntext command for theassociated footnote;
%% use the corref command within \author for corresponding author footnotes;
%% use the cortext command for theassociated footnote;
%% use the ead command for the email address,
%% and the form \ead[url] for the home page:
%% \title{Title\tnoteref{label1}}
%% \tnotetext[label1]{}
%% \author{Name\corref{cor1}\fnref{label2}}
%% \ead{email address}
%% \ead[url]{home page}
%% \fntext[label2]{}
%% \cortext[cor1]{}
%% \address{Address\fnref{label3}}
%% \fntext[label3]{}

\title{Random $k$-noncrossing RNA Structures}

%% use optional labels to link authors explicitly to addresses:
%% \author[label1,label2]{}
%% \address[label1]{}
%% \address[label2]{}

\author{William Y.C. Chen , Hillary S.~W. Han and Christian M. Reidys$^{*}$}

\address{$^*$Center for Combinatorics, LPMC-TJKLC\\
         Nankai University  \\
         Tianjin 300071\\
         P.R.~China\\
         Phone: *86-22-2350-6800\\
         Fax:   *86-22-2350-9272\\
         duck@santafe.edu}

\begin{abstract}
%% Text of abstract
In this paper we derive polynomial time algorithms that generate
random $k$-noncrossing matchings and $k$-noncrossing RNA structures
with uniform probability. Our approach employs the bijection between
$k$-noncrossing matchings and oscillating tableaux and the
$P$-recursiveness of the cardinalities of $k$-noncrossing matchings.
The main idea is to consider the tableaux sequences as paths of
stochastic processes over shapes and to derive their transition
probabilities.
\end{abstract}

\begin{keyword}
%% keywords here, in the form: keyword \sep keyword
RNA pseudoknot structure \sep $k$-noncrossing structure \sep uniform
generation \sep oscillating tableaux \sep stochastic process
%% PACS codes here, in the form: \PACS code \sep code

%% MSC codes here, in the form: \MSC code \sep code
%% or \MSC[2008] code \sep code (2000 is the default)

\end{keyword}

\end{frontmatter}

%% \linenumbers

%% main text

%%%
%%%%%%%%%%%%%%%%%%%%%%%%%%%%%%%%%%%%%%%%%%%%%%%%%%%%%%%%%%%%%%%%%%%%%%%%
%%%
\section{Introduction}
In this paper we generate random $k$-noncrossing partial matchings
and $k$-noncrossing RNA structures with uniform probability in
polynomial time.

Three decades ago Waterman pioneered the concept of RNA secondary
structures \cite{Waterman79,Waterman94}. These coarse grained RNA
structures are subject to strict combinatorial constraints: there
exist no two arcs that cross in their diagram representation. It is
well-known, however, that there exist cross-serial interactions,
i.e.~crossing base pairs \cite{Science:05a}. These configurations
are called pseudoknots \cite{Westhof:92a} and occur in functional
RNA (RNAseP \cite{Loria:96a}), ribosomal RNA \cite{Konings:95a} and
are conserved in the catalytic core of group I introns. Pseudoknots
appear in plant viral RNAs and {\it in vitro} RNA evolution
\cite{Tuerk:92} experiments have produced families of RNA structures
with pseudoknot motifs, when binding HIV-1 reverse transcriptase.
$k$-noncrossing RNA structures \cite{Reidys:08lego} allow to express
pseudoknots and generalize the concept of RNA secondary structures
in a natural way. Due to their cross-serial interactions pseudoknot
structures cannot be recursively generated. This renders the {\it ab
initio} folding into minimum free energy configurations
\cite{Fenix:08} as well as the derivation of detailed statistical
properties a difficult task.

A partial matching over $[n]=\{1,\dots,n\}$ is a labeled graph on
$[n]$, having vertices of degree at most one, represented by drawing
its vertices in increasing order on a horizontal line and its arcs
$(i,j)$ in the upper halfplane. Without loss of generality we shall
assume $i<j$. Two arcs $(i_1,j_1)$ and $(i_2,j_2)$ are crossing if
$i_1<i_2<j_1<j_2$ holds and nesting if $i_1<i_2<j_2<j_1$. A
$k$-crossing is a sequence of arcs $(i_1,j_1),\dots,(i_k,j_k)$ such
that $i_1<i_2<\dots<i_k<j_1<j_2<\dots <j_k$. There is an analogous
notion of a $k$-nesting. A partial matching is called
$k$-noncrossing ($k$-nonnesting) if it does not contain any
$k$-crossing ($k$-nesting). A partial matching without isolated
vertices is called a matching.
%%%
%%%%%%%%%%%%%%%%%%%%%%%%%%%%%%%%%%%%%%%%%%%%%%%%%%%%%%%%%%%%%%%%%%%%%%%
%%%
\begin{figure}[ht]
\centerline{\epsfig{file=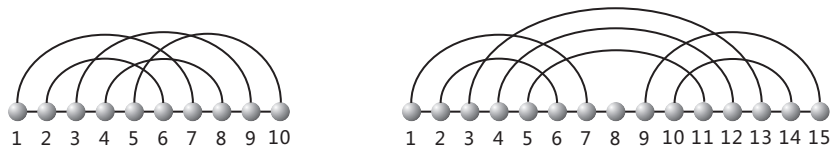,width=0.8\textwidth}\hskip15pt}
\caption{\small $k$-noncrossing diagrams: a $4$-noncrossing (left)
and a $3$-noncrossing diagram (right).} \label{F:dia}
\end{figure}
%%%
%%%%%%%%%%%%%%%%%%%%%%%%%%%%%%%%%%%%%%%%%%%%%%%%%%%%%%%%%%%%%%%%%%%%%%
%%%
The numbers of $k$-noncrossing partial matchings and $k$-noncrossing
matchings over $[n]$ are denoted by $f_k^*(n)$ and $f_k(n)$,
respectively. A $k$-noncrossing RNA structure
\cite{Reidys:07pseu,Reidys:08lego} is a $k$-noncrossing partial
matching without arcs of the form $(i,i+1)$, $1\le i\le n-1$.

\begin{figure}[ht]
\centerline{\epsfig{file=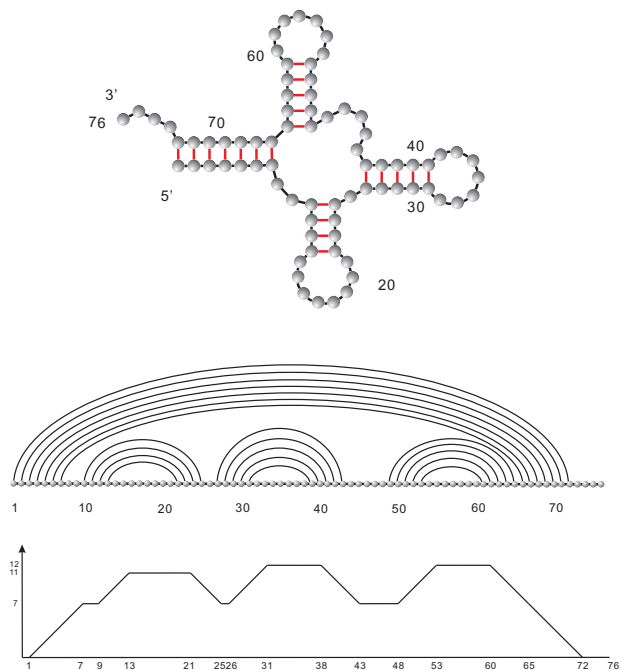,width=0.6\textwidth}\hskip15pt}
\caption{\small The phenylalanine tRNA secondary structure
represented as a planar graph (top), $2$-noncrossing diagram
(middle) and Motzkin-path (bottom).} \label{F:sec}
\end{figure}

For $k$-noncrossing partial matchings there exists no relation
expressing $k$-noncrossing partial matchings over $n$ vertices via
those over $j<n$ vertices as, for instance, the path-concatenation
formula expressing Motzkin-paths, see Figure~\ref{F:sec}, of length
$n$ via shorter paths depending on whether the initial step is an
up-step or a horizontal-step.
%%%
%%%%%%%%%%%%%%%%%%%%%%%%%%%%%%%%%%%%%%%%%%%%%%%%%%%%%%%%%%%%%%%%%%%%%%%
%%%
%%%
%%%%%%%%%%%%%%%%%%%%%%%%%%%%%%%%%%%%%%%%%%%%%%%%%%%%%%%%%%%%%%%%%%%%%%
%%%
Flajolet {\it et al.} \cite{Flajolet} have shown how to uniformly
generate via Boltzman generators elements of a combinatorial class,
for which such a recurrence exists. However, there is no comparable
framework for the uniform generation of elements a non-inductive
combinatorial class, a question arguably of pure mathematical
interest. The subject of uniform generation via Markov-processes has
been studied extensively. A computational study on the uniform
generation of random graphs via Markov-chains has been given by
\cite{Barbosa}. Work on the uniform generation of specific graphs in
the context of parallel random access machine (PRAM) can be found in
\cite{Zito} and Jerrum {\it et al.}
\cite{Vazirant,Jerrum:90,Jerrum:891} studied approximation
algorithms in the context of rapidly mixing Markov-chains
\cite{Aldous:82}.

\begin{figure}[ht]
\centerline{\epsfig{file=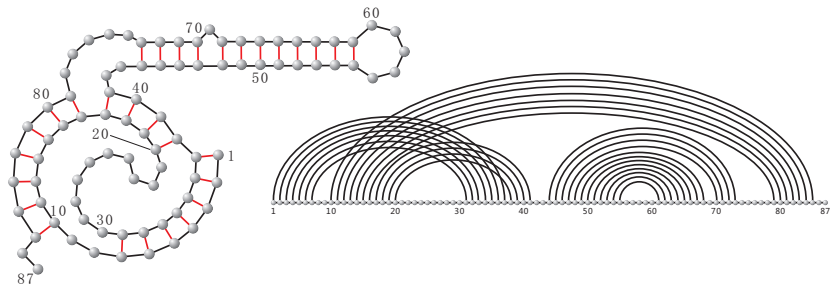,width=0.8\textwidth}\hskip15pt}
\caption{\small The Hepatitis Delta Virus (HDV)-pseudoknot structure
represented as a planar graph and as a diagram: we display the
$3$-noncrossing structure as folded by the {\it ab initio} folding
algorithm {\sf cross} \cite{Fenix:08} (left) and the diagram
representation (right).} \label{F:pse}
\end{figure}

The motivation of this paper comes from the above mentioned RNA
pseudoknot structures \cite{Science:05a}, modeled as $k$-noncrossing
partial matchings without $1$-arcs
\cite{Reidys:07pseu,Reidys:08lego}, see Figure~\ref{F:pse}.
%%%
%%%%%%%%%%%%%%%%%%%%%%%%%%%%%%%%%%%%%%%%%%%%%%%%%%%%%%%%%%%%%%%%%%%%%%
%%%
%%%
%%%%%%%%%%%%%%%%%%%%%%%%%%%%%%%%%%%%%%%%%%%%%%%%%%%%%%%%%%%%%%%%%%%%%%%
%%%
At present time only a few statistical results, that is, central
limit and discrete limit theorems, derived via singularity analysis
of the corresponding bivariate generating functions are known
\cite{Reidys:limit,Fenix,Han,mathhit}. The results of this paper do
not only facilitate the derivation of detailed statistics of RNA
pseudoknot structures but also open the door for studies along the
lines of \cite{7280} and novel, randomized folding routines. Our
algorithms are freely available at {\tt
http://www.combinatorics.cn/cbpc/unif.html}.

Our approach is as follows: we consider the bijection \cite{Chen}
between $k$-noncrossing partial matchings and specific sequences of
Ferrers diagrams and then use the $D$-finiteness \cite{Stanley} of
the ordinary generating function. In some sense, $D$-finiteness is
the next best thing if constructive recurrences are not available.
$D$-finiteness implies $P$-recursiveness \cite{Stanley}, i.e.~the
existence of a finite recurrence relation for the cardinalities of
the combinatorial class with polynomial coefficients. Therefore the
key quantities, i.e.~the transition probabilities of the specific
stochastic processes can be derived with linear time complexity.

%%%
%%%%%%%%%%%%%%%%%%%%%%%%%%%%%%%%%%%%%%%%%%%%%%%%%%%%%%%%%%%%%%%%%%%%%%%%
%%%

\label{S:bijection}\section{The bijection}

%%%
%%%%%%%%%%%%%%%%%%%%%%%%%%%%%%%%%%%%%%%%%%%%%%%%%%%%%%%%%%%%%%%%%%%%%%%%
%%%

In this section we recall the main ideas on the bijection between
$k$-noncrossing partial matchings and $*$-tableaux, a specific class
of vacillating tableaux \cite{Chen}. The bijection facilitates the
interpretation of $k$-noncrossing partial matchings and
$k$-noncrossing structures as paths of the stochastic processes.

A Ferrers diagram, or shape, $\lambda$, is a collection of squares
arranged in left-justified rows with weakly decreasing number of
boxes in each row. A standard Young tableau (SYT), denoted by $T$,
is a filling of the squares by numbers which is strictly decreasing
in each row and in each column. We refer to standard Young tableaux
as Young tableaux. Elements can be inserted into SYT via the
RSK-algorithm \cite{ec1}. We will refer to SYT simply as tableaux. A
$*$-tableaux of shape $\lambda$ and of length $n$ is a sequence of
shapes $(\lambda^i)_{i=0}^n$, where $\lambda^0=\varnothing$,
$\lambda^n=\lambda$, such that for all $ 1\leq i \leq n$, the shape
$\lambda^i$ is obtained from $\lambda^{i-1}$ by either adding one
square, removing one square, or doing nothing. A $*$-tableaux, in
which any two subsequent shapes $\lambda^{i-1},\lambda^i$ are
different is called an oscillating tableaux.

Our first observation \cite{Chen} puts RSK-insertion into context
with $*$-tableaux. It is the key for proving Theorem~\ref{T:chen},
which establishes the bijection between $*$-tableaux of empty shape
and length $n$, having at most $(k-1)$ rows and $k$-noncrossing,
partial matchings on $[n]$. It may be viewed as a ``reverse'' RSK,
facilitating the construction of $*$-tableaux via partial matchings.

\begin{figure}[ht]
\centerline{\epsfig{file=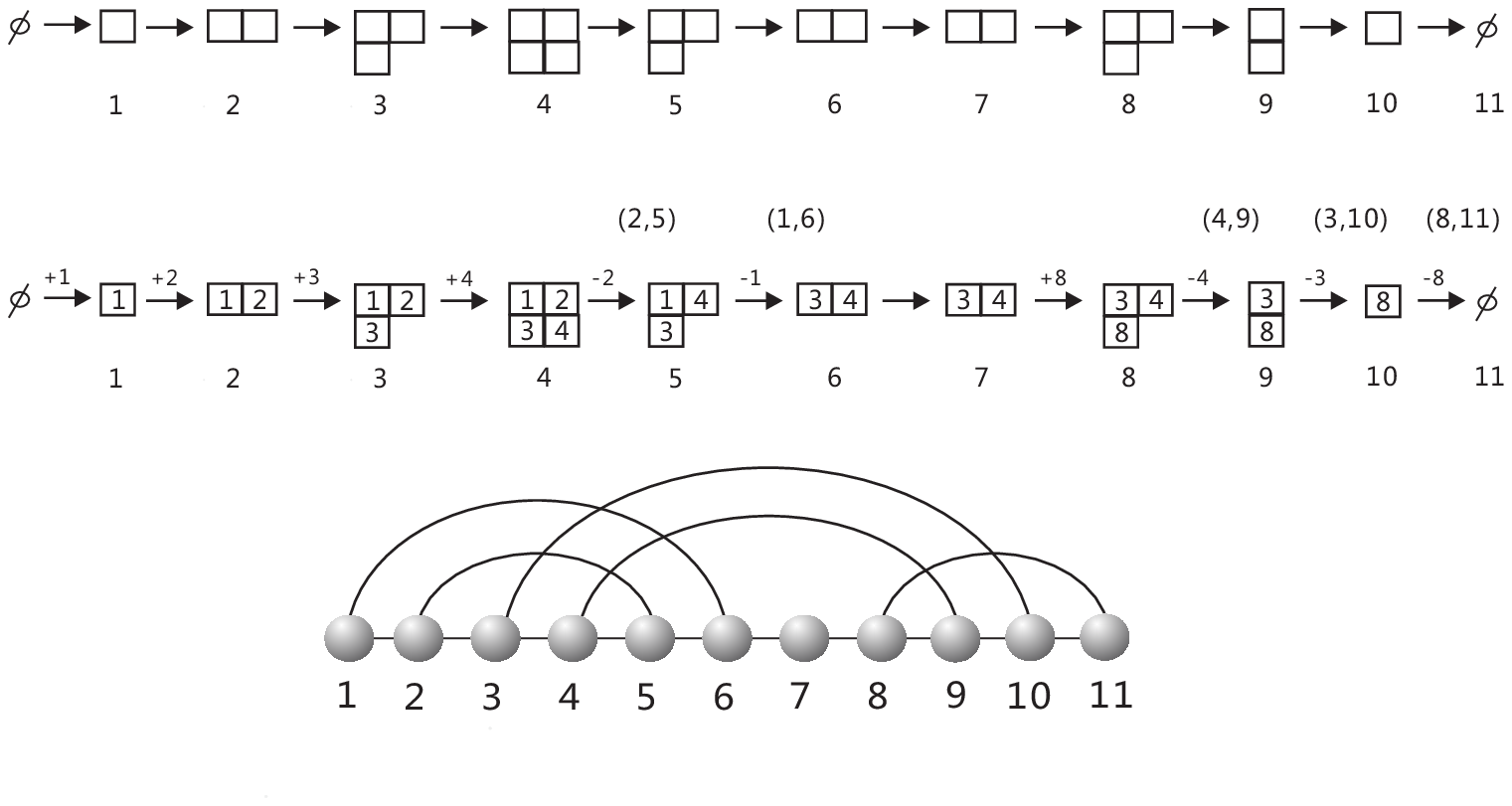,width=0.6\textwidth}\hskip15pt}
\caption{From $*$-tableaux to partial matchings, see
Lemma~\ref{L:extract}. If $\lambda^i\setminus
\lambda^{i-1}=-\square$, then the unique number is extracted, which,
if RSK-inserted into $\lambda^{i},$ recovers $\lambda^{i-1}$. This
yields the arc-set of a $k$-noncrossing, partial
matching.}\label{F:tab-pm}
\end{figure}

%%%
%%%%%%%%%%%%%%%%%%%%%%%%%%%%%%%%%%%%%%%%%%%%%%%%%%%%%%%%%%%%%%%%%%%%%%
%%%
\begin{lemma}\cite{Chen,Reidys:07vac}\label{L:extract}
Suppose we are given two shapes
$\lambda^{i}\subsetneq\lambda^{i-1}$, which differ by exactly one
square. Let $T_{i-1}$ and $T_{i}$ be SYT of shape $\lambda^{i-1}$
and $\lambda^{i}$, respectively. Then there exists a unique $j$
contained in $T_{i-1}$ and a unique tableau $T_{i}$ such that
$T_{i-1}$ is obtained from $T_{i}$ by inserting $j$ via the
$\text{\rm RSK}$-algorithm.
\end{lemma}
%%%
%%%%%%%%%%%%%%%%%%%%%%%%%%%%%%%%%%%%%%%%%%%%%%%%%%%%%%%%%%%%%%%%%%%%%%
%%%
We shall proceed by illustrating how the bijection \cite{Chen}
works. Given a $*$-tableaux of empty shape,
$(\varnothing,\lambda^1,\dots, \lambda^{n-1},\varnothing)$, reading
$\lambda^i\setminus\lambda^{i-1}$
from left to right, at step $i$, we do the following:\\
$\bullet$ for a $+\square$-step we insert $i$ into the new square\\
$\bullet$ for a $\varnothing$-step we do nothing\\
$\bullet$ for a $-\square$-step we extract the unique entry, $j(i)$,
    of the tableaux
    $T^{i-1}$, which via RSK-insertion into $T^i$ recovers it
    (Lemma~\ref{L:extract}).
The latter extractions generate the arc-set $\{(i,j(i))\mid i
\;\text{\rm is a $-\square$-step}\}$ of a $k$-noncrossing, partial
matching, see Fig.~\ref{F:tab-pm}.
%%%
%%%%%%%%%%%%%%%%%%%%%%%%%%%%%%%%%%%%%%%%%%%%%%%%%%%%%%%%%%%%%%%%%%%%%
%%%
%%%%%%%%%%%%%%%%%%%%%%%%%%%%%%%%%%%%%%
%%%%%%%%%%%%%%%%%%%%%%%%%%%%%%%%%%%%%%%%%%%%%%%%%%%%%%%%%%%%%%%
%%%%%%%%%%%%%%%%%%%%%%%%%%%%%%%%5
%%%%%%%%%%%%%%%%%%%%%%%%%%%%%%%%%%%%%%%%%%%%%%%%%%%%%%%%%%%%%%%
%%%%%%%%%%%%%%%%%%%%%%%%%%%%%%%%5
%%%

\begin{figure}[ht]
\centerline{\epsfig{file=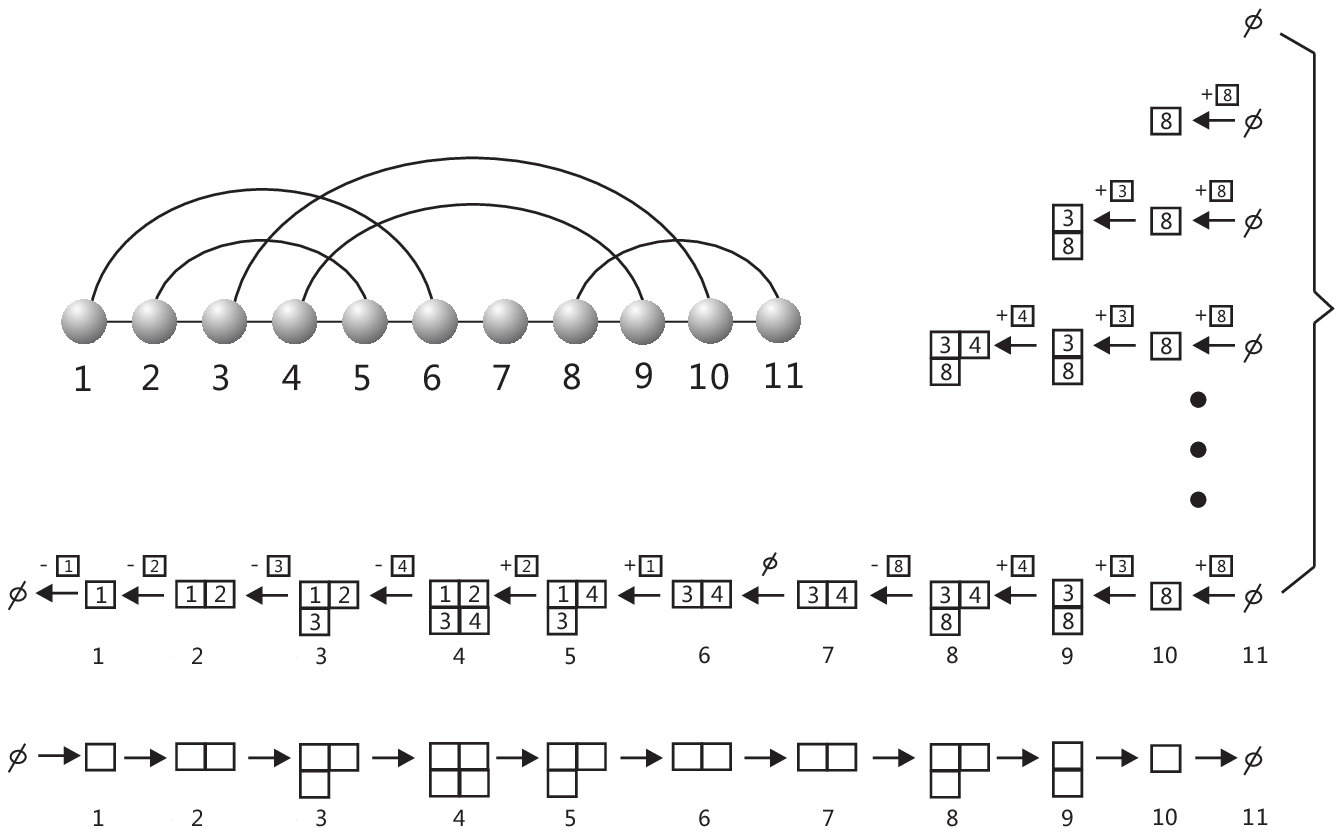,width=.6\textwidth}\hskip15pt}
\caption{From partial matchings to $*$-tableaux via RSK insertion of
the origins of arcs.}\label{F:pm-tab}
\end{figure}

%%%
Given a $k$-noncrossing partial matching, we next construct a unique
$*$-tableaux as follows: starting with the empty shape, consider
the sequence $(n,n-1,\dots,1)$ and do the following:\\
$\bullet$ if $j$ is the endpoint of an arc $(i,j)$, then RSK-insert $i$\\
$\bullet$ if $j$ is the startpoint of an arc $(j, s)$, then remove
the square containing $j$.\\
$\bullet$ if $j$ is an isolated point, then do nothing,
see Fig.~\ref{F:pm-tab}. \\
%%%
%%%%%%%%%%%%%%%%%%%%%%%%%%%%%%
%%%%%%%%%%%%%%%%%%%%%%%%%%%%%%%%%%%%%%%%%%%%%%%%%%%%%%%%%%%%%%%%%%%%%
%%%
%%%%%%%%%%%%%%%%%%%%%%%%%%%%%%%%%%%%%%%%%%%%%%%%%%%%%%%%%%%%%%%%%%%%%

%%%
The above construction leads to
\begin{theorem}\cite{Chen}\label{T:chen}
Each $*$-tableaux of length $n$, containing shapes with at most
$(k-1)$-rows, corresponds uniquely to a $k$-noncrossing partial
matching on $[n]$.
\end{theorem}
%%%
%%%%%%%%%%%%%%%%%%%%%%%%%%%%%%%%%%%%%%%%%%%%%%%%%%%%%%%%%%%%%%%%%%%%%
%%%

Of course, the above bijection induces a correspondence between
oscillating tableaux and $k$-noncrossing matchings.

%%%
%%%%%%%%%%%%%%%%%%%%%%%%%%%%%%%%%%%%%%%%%%%%%%%%%%%%%%%%%%%%%%%%%%%%%
%%%

\label{S:Weyl}\section{D-finiteness}

%%%
%%%%%%%%%%%%%%%%%%%%%%%%%%%%%%%%%%%%%%%%%%%%%%%%%%%%%%%%%%%%%%%%%%%%%
%%%
Suppose $(\lambda^i)_{i=0}^n$ is a $*$-tableaux of shape $\lambda$
having at most $(k-1)$ rows. Let ${\sf O}_k^*(\lambda^i,n-i)$ and
${\sf O}^0_k (\lambda^i,n-i)$ denote the numbers of $*$-tableaux and
oscillating tableaux of shape $\lambda^i$ and length $(n-i)$,
respectively. In this section we establish that these quantities can
be computed with $O(n)$ time and space complexity. In addition, in
the case of $k=3$, we derive explicit formulas.

%%
%%%%%%%%%%%%%%%%%%%%%%%%%%%%%%%%%%%%%%%%%%%%%%%%%%%%%%%%%%%%%%%%%%%%%%%%
%%%
\subsection{The exponential generating function}
%%%
%%%%%%%%%%%%%%%%%%%%%%%%%%%%%%%%%%%%%%%%%%%%%%%%%%%%%%%%%%%%%%%%%%%%%%%%
%%%
Given a $*$-tableaux of shape $\lambda$, $(\lambda^i)_{i=0}^n$,
where $\lambda^n=\lambda$, we consider the number of squares in the
$s$th row of shape $\lambda^i$, denoted by $x_s(i)$. It is evident
that a $*$-tableaux of shape $\lambda$ with at most $(k-1)$ rows
uniquely corresponds to a walk of length $n$ which starts at
$a=(k-1,k-2, \ldots, 1)$ and ends at
$b=(k-1+x_1(n),\dots,1+x_{k-1}(n))$ having steps $0, \pm e_i$, $1
\leq i \leq k-1$ such that $0 < x_{k-1} < \ldots < x_1$ at any step.
That is, a $*$-tableaux of shape $\lambda$ with at most $(k-1)$ rows
corresponds to a lattice path in $\mathbb{Z}^{k-1}$ that remains in
the interior of the dominant Weyl chamber \cite{Zeilberger}. For
$a,b \in \mathbb{Z}^{k-1}$, let $\gamma^*(a,b)$ denote a walk of
length $n$ which starts at $a=(k-1,k-2, \ldots, 1)$, ends at $b$ and
that has steps $0, \pm e_i$, $1\leq i \leq k-1$ such that $0 <
x_{k-1} < \ldots < x_1$ at any step. Let $\Gamma_{n}^*(a,b)$ be the
number of such walks. For our purposes, it suffices to consider
walks without $0$-steps. To this end, let $\gamma_0^*(a,b)$ denote a
$\gamma^*(a,b)$-walk that does not contain any zero-steps and let
$\Gamma_{n}^0(a,b)$ denote their number. In case of $a=b=(k-1,
\ldots, 1)$, Theorem~\ref{T:chen} implies
\begin{equation}
\Gamma_{n}^{*}(a,b)={\sf O}_k^*(\lambda,n) \quad \text{\rm and}\quad
{\Gamma}_{n}^0(a,b)={\sf O}_k^0(\lambda,n),
\end{equation}
where $\lambda$ represents the unique shape with at most $(k-1)$
rows that corresponds to the lattice point $b\in\mathbb{Z}^{k-1}$.
Let $I_r(2x)=\sum_{j \geq 0} \frac{x^{2j+r}}{j! (r+j)!}$ be the
hyperbolic Bessel function of the first kind of order $r$. The
following lemma is implied by the reflection principle \cite{Andre}
and due to Grabiner and Magyar \cite{Grabiner}. It expresses the
exponential generation functions of $\Gamma_{n}^{*}(a,b)$ and
${\Gamma}_{n}^0(a,b)$ via a determinant of Bessel functions:
%%%
%%%%%%%%%%%%%%%%%%%%%%%%%%%%%%%%%%%%%%%%%%%%%%%%%%%%%%%%%%%%%%%%%%%%%%%%%%%%
%%%
\begin{lemma}\cite{Grabiner}\label{L:Gamma}
The exponential generating function for the numbers of
$k$-noncrossing matchings, $\Gamma_n^{0}(a,b)$, is given by
\begin{equation}\label{E:G+}
\sum_{n \geq 0} \Gamma_n^{0}(a,b)\frac{x^n}{n!}= \det[
I_{a_i-b_j}(2x)-I_{a_i+b_j}(2x) ]|_{i,j=1}^{k-1}.
\end{equation}
\end{lemma}
%%%
%%%%%%%%%%%%%%%%%%%%%%%%%%%%%%%%%%%%%%%%%%%%%%%%%%%%%%%%%%%%%%%%%%%%%%%%%%%%
%%%
Consequently we have an algebraic relation for the exponential
generating function. Since $D$-finite functions form an algebra
\cite{Stanley}, the above relation implies, that the ordinary
generating function $\sum_{n \geq 0} \Gamma_n^{0}(a,b)x^n$ is also
$D$-finite. That is, the ordinary generating function of oscillating
tableaux with at most $(k-1)$ rows of arbitrary shape $\lambda$,
$\sum_{n \geq 0}{\sf O}_k^0(\lambda,n) x^n$, is $D$-finite
\cite{Stanley}. Since $D$-finitness is equivalent to the
$P$-recursiveness \cite{ec1}, we have
%%%
%%%%%%%%%%%%%%%%%%%%%%%%%%%%%%%%%%%%%%%%%%%%%%%%%%%%%%%%%%%%%%%%%%%%%%%%%%
%%%
\begin{corollary}\label{C:recursive}
For fixed shape $\lambda$ with at most $(k-1)$ rows and
$n\in\mathbb{N}$, there exists some $m\in\mathbb{N}$ and polynomials
$p_0(n),\dots,p_m(n)$ such that
\begin{equation}\label{E:recur}
p_m(n){\sf O}_k^0(\lambda,n+m)+\dots +p_0(n){\sf O}_k^0(\lambda,n)
=0.
\end{equation}
In particular, the numbers ${\sf O}_k^0(\lambda,n)$ can be computed
in $O(n)$ time.
\end{corollary}
%%%
%%%%%%%%%%%%%%%%%%%%%%%%%%%%%%%%%%%%%%%%%%%%%%%%%%%%%%%%%%%%%%%%%%%%%%%%%%
%%%
The key point here is, that for given $n$ and $\lambda$, the
derivation of eq.~(\ref{E:recur}) is a preprocessing step. It has to
be derived only once, for instance employing Zeilberger's algorithm
\cite{Zeilberger-algo, Salvy-Zimmerman}.

%%%
%%%%%%%%%%%%%%%%%%%%%%%%%%%%%%%%%%%%%%%%%%%%%%%%%%%%%%%%%%%%%%%%%%%%%%%%
%%%
\subsection{The case of 3-noncrossing partial matchings}
$\,$ Let $(\pi_1,\pi_2)$ be a $\mathscr{D}$-pair with the endpoints
$V_u=(n,h_u)$, $u=1,2$ and let $t(n,h_1,h_2)$ denote the number of
such pairs. A $\mathscr{D}$-pair can readily be identified with a
pair of non-intersecting paths, $(\pi_1',\pi_2')$ as follows:
\begin{eqnarray*}
\pi_1' & = & ((i,\pi_1(i)+2))_i \\
\pi_2' & = & \pi_2.
\end{eqnarray*}
Since $(\pi_1',\pi_2')$ are non-intersecting paths, Lindstr\"{o}m's
theorem \cite{Lindstroem,Gessel} allows to compute $t(n,h_1,h_2)$.
Since $\pi_1'(0)=2$, $\pi_2'(0)=0$, $\pi_1'(n)=h_1+2$ and
$\pi_2'(n)=h_2$
\begin{equation}
t(n,h_1,h_2)= p_{(0,2)}^{(n,h_1+2)} p_{(0,0)}^{(n,h_2)} -
p_{(0,2)}^{(n,h_2)} p_{(0,0)}^{(n,h_1+2)}.
\end{equation}
Inserting two up-steps, we observe that $\pi'_1$ uniquely
corresponds to a path starting at $(-2,0)$ and ending at $(n,h_1+2)$
of length $n+2$, that does not cross the line $y=0$ and has only up-
and down-steps. Let $F(n,h)$ denote the number of paths of length
$n$ starting at $(0,0)$ and ending at $(n,h)$, having up- and
down-steps and that stay within the first quadrant. Then
\cite{Andre} %%
%%%%%%%%%%%%%%%%%%%%%%%%%%%%%%%%%%%%%%%%%%%%%%%%%%%%%%%%%%%%%%%%%%%%%%%%
%%%
\begin{lemma}\label{L:22}
\begin{equation}
F(n,h)= { n \choose \frac{n-h}{2}}- { n \choose \frac{n-h-2}{2}}.
\end{equation}
\end{lemma}
%%%
%%%%%%%%%%%%%%%%%%%%%%%%%%%%%%%%%%%%%%%%%%%%%%%%%%%%%%%%%%%%%%%%%%%%%%%%
%%%

{\it Proof:}
 Clearly, there are $\binom{n}{\frac{n-h}{2}}$ paths having
up- and down-steps, that start at $(0,0)$ and end at $(n,h)$. We
call such a path good, if it never touches the line $y=-1$ and bad,
otherwise. Reflecting the segment of a bad path which starts at
$(0,0)$ and ends at the first intersection point at the line $y=-1$,
we observe that the set of bad paths is in one-to-one correspondence
with the set of all paths having only up- and down-steps from
$(0,-2)$ to $(n,h)$. Subtracting the number of these paths,
$\binom{n}{\frac{n-h-2}{2}}$, from the number of all paths, the
lemma follows.

%%%%%%%%%%%%%%%%%%%%%%%%%%%%%%%%%%%%%
We can now give explicit formulas for ${\sf O}^*_3(\lambda^{i},
n-i)$ and ${\sf O}^0_3(\lambda^{i}, n-i)$ \cite{Gouyou,Gessel}

%%%
%%%%%%%%%%%%%%%%%%%%%%%%%%%%%%%%%%%%%%%%%%%%%%%%%%%%%%%%%%%%%%%%%%%%%%%%
%%%
%%%
%%%%%%%%%%%%%%%%%%%%%%%%%%%%%%%%%%%%%%%%%%%%%%%%%%%%%%%%%%%%%%%%%%%%%%%%
%%%
%%%%%%%%%%%%%%%%%%%%%%%%%%%%%%%%%%%%%%%%%%%%%%%%%%%%%%%%%%%%%%%%%%%%%%%%
%%%
\begin{corollary}\label{C:1}
Let $\lambda_{h_1,h_2}^i$ denote the shape with at most two rows,
where $x_1^{\lambda^i}(n)+x_2^{\lambda^i}(n)=h_1$ and
$x_1^{\lambda^i}(n)- x_2^{\lambda^i}(n)=h_2$. Then we have
\begin{equation}\label{E:X0}
\begin{split}
{\sf O}_3^0(\lambda_{h_1,h_2}^i,n-i) & = t(n-i, h_1, h_2) \\
& = F((n-i)+2,h_1+2) F((n-i),h_2) -\\
& \quad \ F((n-i)+2, h_2) F(n-i, h_1+2)
\end{split}
\end{equation}
\begin{equation}\label{E:X1}
\begin{split}
{\sf O}^*_3(\lambda^{i}_{h_1, h_2}, n-i)  =
\begin{cases}
\sum_{l=0}^{\frac{n}{2}}{n-i \choose 2l} t(n-i-2l, h_1, h_2), \\
\hspace*{\fill} \text{\rm for
$(n-i)$ even}\\
\sum_{l=0}^{\frac{n}{2}}{n-i \choose 2l+1} t(n-i-2l-1, h_1, h_2),\\
\hspace*{\fill} \text{\rm for $(n-i)$ odd.}
\end{cases}
\end{split}
\end{equation}
\end{corollary}
%%%
%%%%%%%%%%%%%%%%%%%%%%%%%%%%%%%%%%%%%%%%%%%%%%%%%%%%%%%%%%%%%%%%%%%%%%%%
%%%

%%%
%%%%%%%%%%%%%%%%%%%%%%%%%%%%%%%%%%%%%%%%%%%%%%%%%%%%%%%%%%%%%%%%%%%%%%%%
%%%
%%
%%%%%%%%%%%%%%%%%%%%%%%%%%%%%%%%%%%%%%%%%%%%%%%%%%%%%%%%%%%%%%%%%%%%%%%%%%%%
%%%

\section{Random k-noncrossing partial matchings}\label{S:a1}

%%%
%%%%%%%%%%%%%%%%%%%%%%%%%%%%%%%%%%%%%%%%%%%%%%%%%%%%%%%%%%%%%%%%%%%%%%%%%%%%
%%%
%%%
In this section we generate $k$-noncrossing partial matchings with
uniform probability. The construction is as follows: first we
compute for any shape $\lambda$, having at most $(k-1)$ rows, the
recursion relation of Corollary~\ref{C:recursive}. Second we compute
the array $({\sf O}^*_k(\lambda^{i},n-i))_{\lambda,(n-i)}$, indexed
by $\lambda$ and $(n-i)$. Then we specify a Markov-process that
constructs a $k$-noncrossing partial matching with uniform
probability with linear time and space complexity.

%%%
%%%%%%%%%%%%%%%%%%%%%%%%%%%%%%%%%%%%%%%%%%%%%%%%%%%%%%%%%%%%%%%%%%%%%%%%
%%%
%%%
%%%%%%%%%%%%%%%%%%%%%%%%%%%%%%%%%%%%%%%%%%%%%%%%%%%%%%%%%%%%%%%%%%%%%%%%
%%%
\begin{theorem}\label{T:1}
Random $k$-noncrossing partial matchings can be generated with
uniform probability in polynomial time. The algorithmic
implementation, see Algorithm \ref{A:1}, has $O(n^{k+1})$
preprocessing time and $O(n^{k})$ space complexity. Each
$k$-noncrossing partial matching is generated with $O(n)$ time and
space complexity.
\end{theorem}
%%%
%%
%%%%%%%%%%%%%%%%%%%%%
%%%%%%%%%%%%%%%%%%%%%%%%%%%%%%%%%%%%%%%%%%%%%%%%%%%%%%%%%%%%%%%%%%%%%%%%
%%%
%%%%%%%%%%%%%%%%%%%%%%%%%%%%%%%%%%%%%%%%%%%%%%%%%%%%%%%%%%%%%%%%%%%%%%%%
%%
\hrule
\begin{algorithm}\label{A:1}{{\small }}
\ \hrule \

\hspace{1mm}{\small $1:$}  $Pascal \leftarrow$  {\sf Binomial}($n$)
(computation of all binomial

\hspace{4mm} coefficients, $B(n,h)$.)

\hspace{1mm}{\small $2:$} ${\it PShape} \leftarrow$ {\sf
ArrayP}(n,k) (computation of ${\sf O}_k^*(\lambda^i,n-i)$, $i=$

\hspace{4mm} \, $0,1, \ldots, n-1$, $\lambda^i$, stored in the $k
\times n$ array,

\hspace{4mm} ${\it PShape}$)

\hspace{1mm}{\small $3:$} {\bf while} {$i < n$} {\bf do}

\hspace{1mm}{\small $4:$} \hspace{2mm} {\bf for} {j from $0$ to
$k-1$} {\bf do}

\hspace{1mm}{\small $5:$} \hspace{4mm} X[j]$\leftarrow$ ${\sf
O}_k^*(\lambda^{i+1}, n-(i+1))$

\hspace{1mm}{\small $6:$} \hspace{4mm} {\it sum} $\leftarrow$ {\it
sum}+X[j]

\hspace{1mm}{\small $7:$} \hspace{2mm} {\bf end for}

\hspace{1mm}{\small $8:$} \hspace{2mm} {\it Shape} $\leftarrow$ {\sf
Random}({\it sum}) ({\sf Random} generates the random

\hspace{1mm}{\small $9:$} \hspace{14mm} shape $\lambda^{i+1}$)

{\small $10:$} \hspace{2mm} $i \leftarrow i+1$

{\small $11:$} \hspace{2mm} Insert {\it Shape} into {\it Tableaux}
(the sequence of shapes).

{\small $12:$} {\bf end while}

{\small $13:$} {\sf Map}({\it Tableaux}) (maps {\it Tableaux} into
its corresponding

\hspace{5mm} partial matching)

\hrule
\end{algorithm}

%%%
%%%
%%%%%%%%%%%%%%%%%%%%%%%%%%%%%%%%%%%%%%%%%%%%%%%%%%%%%%%%%%%%%%%%%%%%%%%%
%%%

\begin{figure}[ht]
\centerline{\epsfig{file=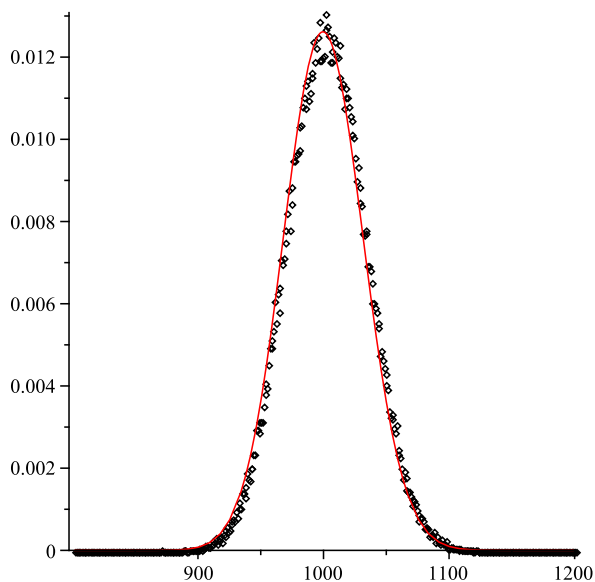,width=0.6\textwidth}\hskip15pt}
\caption{\small Uniformity: for $n=12$ we have $m=f_3^*(12)=99991$
distinct $3$-noncrossing partial matchings. We generate via {\bf
Algorithm 1} $N=10^8$ randomly and display the frequency
distribution of their multiplicities (black dots) and
$\binom{N}{\ell} (1/m)^\ell (1-1/m)^{N-\ell}$ (red curve).}
\label{F:79}
\end{figure}

Figure~\ref{F:79} illustrates that {\bf Algorithm 1} indeed
generates each  $k$-noncrossing partial matching with uniform
probability.
%%%
%%%%%%%%%%%%%%%%%%%%%%%%%%%%%%%%%%%%%%%%%%%%%%%%%%%%%%%%%%%%%%%%%%%%%%%
%%%
%%%
%%%%%%%%%%%%%%%%%%%%%%%%%%%%%%%%%%%%%%%%%%%%%%%%%%%%%%%%%%%%%%%%%%%%%%
%%%

\ {\it Proof:} Suppose $(\lambda^i)_{i=0}^n$ is a $*$-tableaux of
shape $\lambda$ having at most $(k-1)$ rows. By definition, a shape
$\lambda^{i+1}$ does only depend on its predecessor, $\lambda^{i}$.
Accordingly, we can interpret any given $*$-tableaux of shape
$\lambda$ as a path of a Markov-process $(X^i)_{i=0}^n$ over shapes,
given as follows:  \\
$\bullet$ $X^0=X^n=\varnothing$ and $X^i$ is a shape having at most
$(k-1)$
          rows\\
$\bullet$ for $0\le i\le n-1$, $X^i$ and $X^{i+1}$ differ by at most
          one square\\
$\bullet$ the transition probabilities are given by
\begin{equation}\label{E:trans1}
\mathbb{P}(X^{i+1}= \lambda^{i+1} \mid X^i=\lambda^i)=\frac{ {\sf
O}^*_k(\lambda^{i+1},n-(i+1))}{{\sf O}^*_k(\lambda^{i},n-i)}.
\end{equation}

We next observe
\begin{equation}\label{E:oha}
\prod_{i=0}^{n}\mathbb{P}(X^{i+1}=\lambda^{i+1}\mid X^i=\lambda^i)=
\frac{1}{{\sf O}^*_k(\varnothing,n)}=\frac{1}{f_k^*(n)},
\end{equation}
where
\begin{equation}\label{E:X11}
\begin{split}
{\sf O}^*_k(\lambda^{i}, n-i)  =
\begin{cases}
\sum_{l=0}^{\frac{n}{2}}{n-i \choose 2l} {\sf
O}_k^0(\lambda^i,n-i-2l), \\
\hspace*{\fill} \text{\rm for $(n-i)$ even} \\
\sum_{l=0}^{\frac{n}{2}}{n-i \choose 2l+1} {\sf
O}_k^0(\lambda^i,n-i-2l-1), \\
\hspace*{\fill} \text{\rm for $(n-i)$ odd.}
\end{cases}
\end{split}
\end{equation}
Accordingly, the Markov-process, $(X^i)_{i=0}^n$, generates
$k$-noncrossing partial matchings with uniform probability. Clearly,
the Pascal triangle of binomial coefficients can be generated in
$O(n^2)$ time and space and for any fixed $\lambda^i$, having at
most $(k-1)$ rows, we can via Corollary~\ref{C:recursive} compute
${\sf O}_k^0(\lambda^i, n-i)$ in $O(n)$ time. Consequently, we can
generate the array of numbers ${\sf O}_k^0(\lambda^i,n-i)$ as well
as ${\sf O}_k^*(\lambda^i,n-i)$ for all shapes $\lambda$ in
$O(n^2)+O(n)\,O(n) \, O(n^{k-1})$ time and $O(n^{k})$ space. The
first factor $O(n)$ represents the time complexity for deriving the
recursion and the second comes from the computation of all numbers
${\sf O}_k^0(\lambda^i,n-i)$ for fixed $\lambda=\lambda^i$ for all
$(n-i)$.

As for the generation of a random $k$-noncrossing partial matching,
for each shape $\lambda^i$, the transition probabilities can the be
derived in $O(1)$ time. Therefore, a $k$-noncrossing partial
matching can be computed with $O(n)$ time and space complexity,
whence the theorem.\hfill $\square$

%%%%%%%%%%%%%%%%%%%%%%%%%%%%%%%%%%%%%%%%%%%%%%%%%%%%%%%%%%%%%%%%%%%%%%%%%%%%%%
%%%

\section{Random k-noncrossing RNA structures}\label{S:a2}

%%%
%%%%%%%%%%%%%%%%%%%%%%%%%%%%%%%%%%%%%%%%%%%%%%%%%%%%%%%%%%%%%%%%%%%%%%%%%%%%%%
%%%
In this section we generate $k$-noncrossing structures with uniform
probability. The approach is analogous to that of the previous
section, however the stochastic process required here is not a
Markov-process any more. In order to avoid generating arcs of the
form $(i,i+1)$ ($1$-arcs), some kind of one-step memory is needed.

To formalize this intuition, we shall begin by giving the formula
for the number of $k$-noncrossing structures, or equivalently
partial matchings without $1$-arcs \cite{Reidys:07pseu}
\begin{equation}\label{E:S(n)}
 S_k(n)= \sum_{b=0}^{\frac{n}{2}} (-1)^b {n-b \choose b}
f_k^*(n-2b).
\end{equation}
We observe that a $1$-arc corresponds via Theorem~\ref{T:chen} to a
subsequence of shapes $(\lambda^i$, $\lambda^{i+1}$,
$\lambda^{i+2}=\lambda^i)$, obtained by first adding and then
removing a square in the first row. This sequence corresponds to a
pair of steps $(+\square_1,-\square_1)$, where $+\square_1$ and
$-\square_1$ indicate that a square is added and subtracted in the
first row, respectively.

Let ${\mathcal Q}^{*}_k (\lambda^i, n-i, j)$ denote the set of
$*$-tableaux of shape $\lambda^{i}$ of length $(n-i)$ having at most
$(k-1)$ rows containing exactly $j$ pairs $(+\square_1,-\square_1)$
and set ${\sf Q}_k^{*}(\lambda^i, n-i, j)= |{\mathcal
Q}_k^{*}(\lambda^i, n-i, j)|$. Furthermore, let ${\sf
W}_k^*(\lambda^i,n-i)$ denote the number of $*$-tableaux of shape
$\lambda^i$ with at most $(k-1)$ rows of length $(n-i)$ that do not
contain any such pair of steps.

In terms of $*$-tableaux having at most $(k-1)$ rows,
eq.~(\ref{E:S(n)}) can be rewritten as follows $ {\sf
W}_k^{*}(\varnothing,n)= \sum_{b=0}^{\frac{n}{2}} (-1)^b {n-b
\choose b} {\sf O}^*_k(\varnothing,n-2b)$. We proceed by
generalizing this relation from the empty shape to arbitrary shapes.

\begin{figure}[ht]
\centerline{\epsfig{file=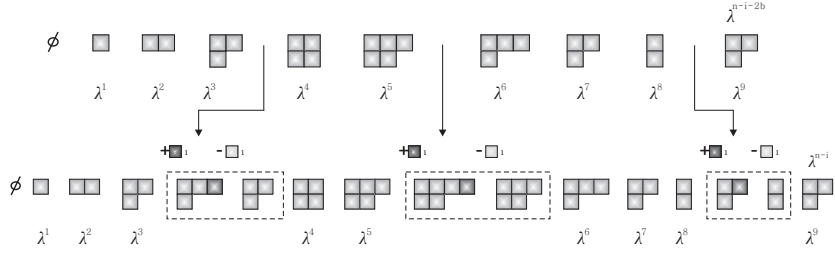,width=.8\textwidth}\hskip15pt}
\caption{Illustration of the proof idea: pairs
$(+\square_1,-\square_1)$ are inserted at positions $3,5$ and $8$,
respectively.}\label{F:square}
\end{figure}

%%%
%%%%%%%%%%%%%%%%%%%%%%%%%%%%%%%%%%%%%%%%%%%%%%%%%%%%%%%%%%%%%%%%%%%%%%%%
%%%
\begin{lemma}\label{L:general}
Let $\lambda^i$ be an arbitrary shape with at most $(k-1)$ rows,
then
\begin{equation}
{\sf W}_k^*(\lambda^i_{},n-i)= \sum_{b=0}^{\frac{n-i}{2}} (-1)^b
{(n-i)-b \choose b} {\sf O}_k^*(\lambda^i,n-i-2b).
\end{equation}
\end{lemma}
%%%
%%%%%%%%%%%%%%%%%%%%%%%%%%%%%%%%%%%%%%%%%%%%%%%%%%%%%%%%%%%%%%%%%%%%%%%%
%%%
{\it Proof:} Let $(\lambda^s)_{s=0}^{(n-2b)-i}$ be a $*$-tableaux of
shape $\lambda^i$. We select from the set $\{0,\dots,(n-2b)-i-1\}$
an increasing sequence of labels $(r_1,\dots,r_b)$. For each $r_s$
we insert a pair $(+\square_1,-\square_1)$ after the corresponding
shape $\lambda^{r_s}$, see Fig.~\ref{F:square}. This insertion
generates a $*$-tableaux
of length $(n-i)$ of shape $\lambda^i$.\\
%%%
%%%%%%%%%%%%%%%%%%%%%%%%%%%%%%%%%%%%
%%%%%%%%%%%%%%%%%%%%%%%%%%%%%%%%%%%%%%%%%%%%%%%%%%%%%%%%%%%%%%%%%%%%%%%%
%%%%%%%%%%%%%%%%%%%%%%
%%%
Considering the above insertion for all sequences $(r_1,\dots,r_b)$,
we arrive at a family $\mathcal{F}_b$ of $*$-tableaux of length
$(n-i)$ containing at least $b$ pairs, $(+\square_1,-\square_1)$.
Since we can insert at any position $0\le h\le  ((n-i)-2b-1)$,
$\mathcal{F}_b$ has cardinality ${(n-i)-b \choose b}{\sf O}_k^{*}
(\lambda^i, n-i-2b)$. By construction, each $*$-tableaux
$(\lambda^s)_{s=0}^{n-i}\in \mathcal{F}_b$, that exhibits exactly
$j$ pairs $(+\square_1,-\square_1)$ appears with multiplicity ${j
\choose b}$, whence

\begin{equation}
\sum_{j \geq b}{j \choose b}{\sf Q}_k^{*}(\lambda^i, n-i,
j)={(n-i)-b \choose b}{\sf O}_k^{*}(\lambda^i, n-i-2b).
\end{equation}
We consider $F_k(x)=\sum_{j \geq 0}{\sf Q}_k^{*}(\lambda^i, n-i, j)
x^j $. Taking the $b$th derivative and setting $x=1$ we obtain $
\frac{1}{b !} \,F_k^{b}(1) =\sum_{j \geq b } {j \choose b}\, {\sf
Q}_k^{*}(\lambda^i, n-i, j)1^{j-b} $ and computing the Taylor
expansion of $F_k(x)$ at $x=1$
\begin{align*}
& F_k(x) =\sum_{b \geq 0}\frac{1}{b !} \,F_k^{b}(1) \,(x-1)^b \\ &=
\sum_{b=0}^{\frac{n-i}{2}}\, {(n-i)-b \choose b}\,{\sf
O}_k^{*}(\lambda^i, n-i-2b)\,(x-1)^b.
\end{align*}

Since ${\sf W}_k^{*}(\lambda^i, n-i)={\sf Q}_k^{*}(\lambda^i, n-i,
0)$ is the constant term of $F_k(x)$, the lemma follows.

We can now prove

%%%
%%%%%%%%%%%%%%%%%%%%%%%%%%%%%%%%%%%%%%%%%%%%%%%%%%%%%%%%%%%%%%%%%%%%%%%%
%%%
\begin{theorem}\label{T:2}
A random $k$-noncrossing structure can be generated, after
polynomial preprocessing time, with uniform probability in linear
time. The algorithmic implementation, see Algorithm~\ref{A:2}, has
$O(n^{k+1})$ pre-processing time and $O(n^{k})$ space complexity.
Each $k$-noncrossing structure is generated with $O(n)$ space and
time complexity.
\end{theorem}
%%%
%%%%%%%%%%%%%%%%%%%%%%%%%%%%%%%%%%%%%%%%%%%%%%%%%%%%%%%%%%%%%%%%%%%%%%%%
%%%
%%%%%%%%%%%%%%%%%%%%%%%%%%%%%%%%%%%%%%%%%%%%%%%%%%%%%%%%%%%%%%%%%%%%%%%%
%%%
\hrule
\begin{algorithm}\label{A:2}{{\small }}
\ \hrule \

\hspace{1mm}{\small $1:$} $Pascal \leftarrow$  {\sf Binomial}($n$)
(computation of all

\hspace{4mm} binomial coefficients, $B(n,h)$.)

\hspace{1mm}{\small $2:$} ${\it PShape} \leftarrow$ {\sf
ArrayP}(n,k) (computation of ${\sf O}_k^*(\lambda^i,n-i), \, i=$

\hspace{4mm} $0,1, \ldots, n-1$, $\lambda^i$)

\hspace{1mm}{\small $3:$} ${\it SShape} \leftarrow$ {\sf
ArrayS}(n,k) (computation of ${\sf W}_k^*(\lambda^i_j, n-i ), \, j=$

\hspace{4mm} $0,1^+,1^-,\dots, (k-1)^+,(k-1)^-; i=0,1, \ldots, n-$

\hspace{4mm} $1$, stored in the $k \times n$ array ${\it SShape}$)

\hspace{1mm}{\small $4:$} {\bf while} {$i< n$} {\bf do}

\hspace{1mm}{\small $5:$} \hspace{2mm} flag $\leftarrow 1$

\hspace{1mm}{\small $6:$} \hspace{2mm} X[0] $\leftarrow$ ${\sf
W}_k^*(\lambda_{0}^{i+1}, n-(i+1))$

\hspace{1mm}{\small $7:$} \hspace{2mm} X[1]$\leftarrow {\sf
W}_k^*(\lambda_{1^+}^{i+1}, n-(i+1))-{\sf
W}_k^*(\lambda_{1^-}^{i+2}, n-(i+2))$

\hspace{1mm}{\small $8:$} \hspace{2mm} {\bf if} {flag=0 and j=2}
{\bf then}

\hspace{1mm}{\small $9:$} \hspace{4mm}  X[2]$\leftarrow 0$

{\small $10:$} \hspace{2mm} {\bf else}

{\small $11:$} \hspace{4mm} X[2]$\leftarrow {\sf
W}_k^*(\lambda_{1^-}^{i+1}, n-(i+1))$

{\small $12:$} \hspace{2mm} {\bf end if}

{\small $13:$} \hspace{2mm} {\it sum} $\leftarrow$ X[0]+X[1]+X[2]

{\small $14:$} \hspace{2mm} {\bf for}{$j$ from $2$ to $k-1$} {\bf
do}

{\small $15:$} \hspace{4mm} X[2j-1] $\leftarrow$ ${\sf
W}_k^*(\lambda_{j^+}^{i+1}, n-(i+1))$

{\small $16:$} \hspace{4mm} X[2j] $\leftarrow$ ${\sf
W}_k^*(\lambda_{j^-}^{i+1}, n-(i+1))$

{\small $17:$} \hspace{4mm} {\it sum}$\leftarrow${\it
sum}+X[2j-1]+X[2j]

{\small $18:$} \hspace{2mm} {\bf end for}

{\small $19:$} \hspace{2mm} {\it Shape} $\leftarrow$ {\sf
Random}({\it sum}) ({\sf Random} generates the random

\hspace{8mm} shape $\lambda^{i+1}_j$ with probability X[j]/{\it
sum})

{\small $20:$} \hspace{2mm} $i \leftarrow i+1$

{\small $21:$} \hspace{2mm} {\bf if} {{\it Shape}
=$\lambda_{1^+}^i$} {\bf then}

{\small $22:$} \hspace{4mm} flag $\leftarrow 0$

{\small $23:$} \hspace{2mm} {\bf end if}

{\small $24:$} \hspace{2mm} Insert $\lambda^{i+1}_j$ into {\it
Tableaux}

{\small $25:$} {\bf end while}

{\small $26:$} {\sf Map}({\it Tableaux}) \ \ \hrule \ \
\end{algorithm}

\begin{figure}[ht]
\centerline{\epsfig{file=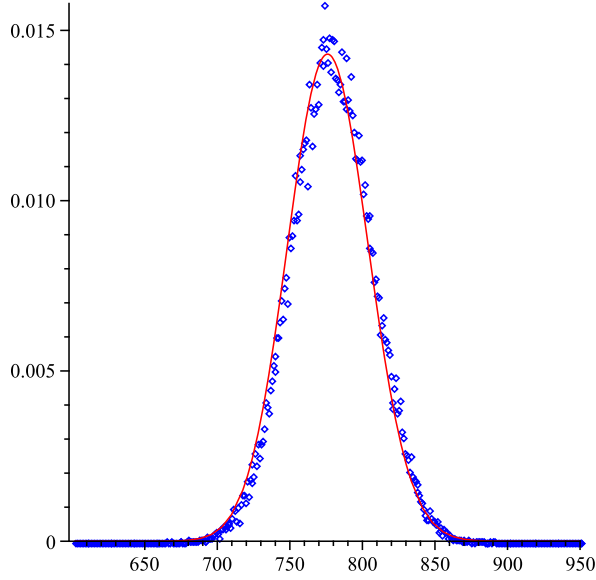,width=0.6\textwidth}\hskip15pt}
\caption{\small Uniformity: for $n=12$ we have $m=S_3(n)=38635$
distinct $3$-noncrossing RNA structures. We generate via {\bf
Algorithm 2} $N=3\times 10^7$ of them and display the frequency
distribution of their multiplicities (blue dots) and
$\binom{N}{\ell} (1/m)^\ell(1-1/m)^{N-\ell}$ (red curve).}
\label{F:80}
\end{figure}

Figure~\ref{F:80} illustrates that {\bf Algorithm 2} generates
$k$-noncrossing RNA structures with uniform probability.
%%%
%%%%%%%%%%%%%%%%%%%%%%%%%%%%%%%%%%%%%%%%%%%%%%%%%%%%%%%%%%%%%%%%%%%%%%%
%%%
%%%
%%%%%%%%%%%%%%%%%%%%%%%%%%%%%%%%%%%%%%%%%%%%%%%%%%%%%%%%%%%%%%%%%%%%%%
%%%

%%%
%%%%%%%%%%%%%%%%%%%%%%%%%%%%%%%%%%%%%%%%%%%%%%%%%%%%%%%%%%%%%%%%%%%%

%%%

{\it Proof:} The idea is to interpret $*$-tableaux without pairs of
steps, $(+\square_1,-\square_1)$, (good $*$-tableaux) as paths of a
stochastic process.

For this purpose we index the shapes $\lambda^{i+1}$ according to
their predecessors: let $i=0,1, \ldots, n-1$ and
$j\in\{0,1^+,1^-,\dots, (k-1)^+,(k-1)^-\}$. Setting
$\lambda^0_j=\varnothing$,
we write $\lambda^{i+1}_{j}$, if $\lambda^{i+1}$ is obtained via \\
$\bullet$ doing nothing ($\lambda_0^{i+1}$)\\
$\bullet$ adding a square in the $j$th row ($\lambda_{j^+}^{i+1}$) \\
$\bullet$ deleting a square in the $j$th row ($\lambda_{j^-}^{i+1}$).\\
With this notation, the number of good $*$-tableaux of shape
$\lambda^{i+1}_{1^+}$ of length $(n-(i+1))$ is given as follows:
\begin{equation*}
{\sf V}_k^*(\lambda^{i+1}_{1^+},n-(i+1))= {\sf
W}_k^*(\lambda^{i+1}_{1^+},n-(i+1))- {\sf
W}_k^*(\lambda^{i+2}_{1^-},n-(i+2)).
\end{equation*}
In order to derive transition probabilities, we establish two
equations: first, for any $\lambda^i_j$, where $j\neq 1^+$
\begin{eqnarray*}
&& {\sf W}_k^*(\lambda^i_j,n-i) \\
& = & {\sf V}_k^*(\lambda^{i+1}_{1^+},n-(i+1))+
{\sf W}_k^*(\lambda^{i+1}_{1^-},n-(i+1))+\\
&& \sum_{h=2}^{k-1}\left({\sf W}_k^*(\lambda^{i+1}_{h^+},n-(i+1))+
{\sf W}_k^*(\lambda^{i+1}_{h^-},n-(i+1))\right)+\\
&& {\sf W}_k^*(\lambda^{i+1}_0,n-(i+1))
\end{eqnarray*}
and second, in case of $j=1^+$
\begin{eqnarray*}
&&{\sf V}_k^*(\lambda^{i}_{1^+},n-i) \\
& = &
{\sf V}_k^*(\lambda^{i+1}_{1^+},n-(i+1))+\\
&& \sum_{h=2}^{k-1}\left({\sf W}_k^*(\lambda^{i+1}_{h^+},n-(i+1))+
{\sf W}_k^*(\lambda^{i+1}_{h^-},n-(i+1))\right)+\\
&& {\sf W}_k^*(\lambda^{i+1}_{0},n-(i+1)).
\end{eqnarray*}
We are now in position to specify the process $(X^i)_{i=0}^n$:\\
$\bullet$ $X^0=X^n=\varnothing$ and $X^i$ is a shape having at most
$(k-1)$
          rows\\
$\bullet$ for $0\le i\le n-1$, $X^i$ and $X^{i+1}$ differ by at most
          one square\\
$\bullet$ there exists no subsequence $X^{i},X^{i+1},X^{i+2}=X^i$
          obtained by first adding and second removing a square in the first
          row\\
$\bullet$ for $j\neq 1^+$
\begin{equation}\label{E:w0}
\mathbb{P}(X^{i+1}  =  \lambda^{i+1}_l \mid X^i=\lambda^i_j)=
\begin{cases}
\frac{{\sf W}^*_k(\lambda^{i+1}_l,n-(i+1))}{{\sf
W}^*_k(\lambda^{i}_j,n-i)}, &
\text{\rm for } l\neq 1^+ \\
\frac{{\sf V}_k^*(\lambda^{i+1}_{1^+},n-(i+1))}{{\sf
W}^*_k(\lambda^{i}_j,n-i)}, & \text{\rm for $l=1^+$}\\
\end{cases}
\end{equation}
$\bullet$ for $j=1^+$
\begin{equation}\label{E:w1}
\mathbb{P}(X^{i+1}= \lambda^{i+1}_l \mid X^i=\lambda^i_{1^+})=
\begin{cases}
\frac{{\sf W}^*_k(\lambda^{i+1}_l,n-(i+1))}
{{\sf V}^*_k(\lambda^{i}_{1^+},n-i)},  \\
\hspace*{\fill}  \text{\rm for $l\neq 1^+, 1^-$} \\
\frac{{\sf V}^*_k(\lambda^{i+1}_{1^+},n-(i+1))}{ {\sf
V}^*_k(\lambda^{i}_{1^+},n-i)}, \quad \text{\rm for $l=1^+$.}
\end{cases}
\end{equation}
As in the proof of Theorem~\ref{T:1} we observe that
eq.~(\ref{E:w0}) and eq.~(\ref{E:w1}) imply
\begin{equation}
\prod_{i=0}^{n-1}\mathbb{P}(X^{i+1}=\lambda^{i+1}\mid
X^i=\lambda^i)= \frac{{\sf W}^*_k(\lambda^{n}=\varnothing,0)} {{\sf
W}^*_k(\lambda^{0}=\varnothing,n)} =\frac{1}{{\sf
W}^*_k(\varnothing,n)}.
\end{equation}
Consequently, the process $(X^i)_{i=0}^n$ generates random
$k$-noncrossing structures with uniform probability in $O(n)$ time
and space. According to Corollary~\ref{C:recursive}, we can for any
$\lambda^i$, having at most $(k-1)$ rows, compute ${\sf
O}_k^0(\lambda^i,n-i)$ in $O(n)$ time. Consequently, we can generate
the arrays $({\sf O}_k^*(\lambda^i,n-i))_{ \lambda^i,n-i}$ and
$({\sf W}^*_k(\lambda^{i},n-i))_{\lambda^i,n-i}$ in
$O(n^2)+O(n^2)\,O(n^{k-1})$ time and $O(n^{k})$ space.

A random $k$-noncrossing structure is then generated as a
$*$-tableaux with at most $(k-1)$ rows using the array $({\sf
W}^*_k(\lambda^{i},n-i))_{\lambda^i,n-i}$ with $O(n)$ time and space
complexity.\hfill $\square$

%\begin{acknowledgments}
This work was supported by the PCSIRT Project of the Ministry of
Education.
%\end{acknowledgments}

\end{document}